\documentclass[11pt]{amsart}
\usepackage{amssymb,xypic,amscd}
\CompileMatrices

\newtheorem{defn}{Definition}
\newtheorem{thm}[defn]{Theorem}

\newtheorem{lem}[defn]{Lemma}
\newtheorem{prop}[defn]{Proposition}

\theoremstyle{remark}
\newtheorem{rem}[defn]{Remark}
\theoremstyle{remark}

\numberwithin{equation}{section} \numberwithin{defn}{section}
\newcommand\glva{{\operatorname{Gl}}(V_A,V_A^+)}

\newcommand\glvae{{\widetilde{\operatorname{Gl}}(V_A,V_A^+)}}
\newcommand\gl{\operatorname{Gl}}

\newcommand\mat{\operatorname{Mat}}

\newcommand\Res{\operatorname{Res}}

\renewcommand\tilde{\widetilde}
\newcommand\ZZ{{\mathbb Z}}
\newcommand\QQ{{\mathbb Q}}
\begin{document}

\title[Steinberg property of the Contou-Carr\`{e}re symbol]{On the Steinberg
property
\\ of the Contou-Carr\`{e}re symbol}
\author{Fernando Pablos Romo}
\address{Departamento de
Matem\'aticas, Universidad de Salamanca, Plaza de la Merced 1-4,
37008 Salamanca, Espa\~na} \email{fpablos@usal.es}
\thanks{
\\ This work is partially supported by the
DGESYC research contract no. BFM2003-00078 and Castilla y Le\'on
regional government contract SA071/04.}

\maketitle

\begin{abstract}  The aim of this work is to show that, when $A$ is an artinian local ring, the Contou-Carr\`{e}re symbol
satisfies the property of Steinberg symbols: $\langle
f,1-f\rangle_{A((t))^\times} = 1$ for all elements $f,1-f \in
A((t))^\times$. Moreover, we offer a cohomological
characterization of the Contou-Carr\`{e}re symbol from the
commutator of a central extension of groups.
\end{abstract}
\maketitle
%%%%%%%%%%%%%%%%%%%%%%%%%%%%%%%%%%%%%%%%%%%%%%%%%%%%%%%%%%%%%
\setcounter{tocdepth}1

%\tableofcontents
%%%%%%%%%%%%%%%%%%%%%%%%%%%%%%%%%%%%%%%%%%%%%%%%%%%%%%%%%%%%%
\section{Introduction}

    In 1971, J. Milnor \cite{Mi} defined the
tame symbol $d_v$ associated with a discrete valuation $v$ on a
field $F$. Explicitly, if $A_v$ is the valuation ring, $p_v$ is
the unique maximal ideal, and $k_v = {A_v}/{p_v}$ is the residue
class field, Milnor defined $d_v \colon F^\times \times F^\times
\to k_v^\times$ by
$$d_v(x,y) = (-1)^{v(x)\cdot v(y)} \frac {x^{v(y)}}{y^{v(x)}}
(\text {mod }p_v).$$

    (Here and below $R^\times$ denotes the multiplicative group of a ring
$R$ with unit).

    J. Milnor proved that the tame symbol is a Steinberg symbol; that is, it is bimultiplicative and
it satisfies the condition that $d_v(x,1-x) = 1$ for all $x\ne 1$.

    In 1994 C. Contou-Carr\`{e}re \cite{Co} defined a
natural transformation greatly generalizing the tame symbol. In
the case of an artinian local base ring $A$ with maximal ideal
$m$, the natural transformation takes the following form. Let
$f,g\in A((t))^\times$ be given, where $t$ is a variable. It is
possible in exactly one way to write:
$$\begin{aligned} f &=
a_0\cdot t^{w(f)} \cdot {\prod}_{i=1}^{\infty} (1-a_{i}t^{i})
\cdot {\prod}_{i=1}^{\infty} (1-a_{-i}t^{-i}) \\ g &= b_0\cdot
t^{w(g)} \cdot {\prod}_{i=1}^{\infty} (1-b_{i}t^{i}) \cdot
{\prod}_{i=1}^{\infty} (1-b_{-i}t^{-i})\,
,\end{aligned}$$\noindent with $w(f),w(g)\in \ZZ$, $a_{
i},b_{i}\in A$ for $i>0$, $a_0, b_0 \in A^\times$, $a_{-i},
b_{-i}\in m$ for $i>0$, and $a_{-i}=b_{-i}=0$ for $i\gg 0$. By
definition, the value of the Contou-Carr\`{e}re symbol is:
$$\langle f,g\rangle_{A((t))^\times}  := (-1)^{w(f)w(g)}\frac
{a_0^{w(g)}{\prod}_{i=1}^{\infty}{\prod}_{j=1}^{\infty}\big (1 -
a_i^{j/(i,j)}b_{-j}^{i/(i,j)}\big
)^{(i,j)}}{b_0^{w(f)}{\prod}_{i=1}^{\infty}{\prod}_{j=1}^{\infty}\big
(1 - a_{-i}^{j/(i,j)}b_{j}^{i/(i,j)}\big )^{(i,j)}} \in
A^\times.$$ The definition makes sense because only finitely many
of the terms appearing in the infinite products differ from 1. The
symbol $\langle\cdot,\cdot\rangle_{A((t))^\times}$ is clearly
antisymmetric and, although is not immediately obvious  from the
definition, it is also bimultiplicative.

     G. W. Anderson and the author
\cite{AP} have interpreted the Contou-Carr\`{e}re symbol $\langle
f,g\rangle_{A((t))^\times}$---up to signs--- as a commutator of
liftings of $f$ and $g$ to a certain central extension of a group
containing $A((t))^\times$, and they have exploited the commutator
interpretation to prove in the style of Tate \cite{Ta} a
reciprocity law for the Contou-Carr\`{e}re symbol on a
non-singular complete curve defined over an algebraically closed
field $k$, $A$ being an artinian local $k$-algebra.

    Moreover, the author has obtained a similar result for an
algebraic curve over a perfect field \cite{Pa2}, and A. Beilinson,
S. Bloch and H. Esnault \cite{BBE} have defined the
Contou-Carr\`{e}re symbol as the commutator pairing in a
Heisenberg super extension. This symbol has also played an
important role in a recent work by M. Kapranov and E. Vasserot
\cite{KV}. In fact, since the Contou-Carr\`{e}re symbol contains
the classical residue, the tame symbol and the Witt residue as
special cases, it is currently an important tool for studying
several topics in Arithmetic Algebraic Geometry.

    An open problem regarding this symbol has been to determine
whether it satisfies the Steinberg property:
\begin{equation} \label{eq:Steinberg} \langle f,1-f\rangle_{A((t))^\times} = 1
\text{ for all } f,1-f \in A((t))^\times\, .\end{equation}

    Here we solve this problem and we show that the Contou-Carr\`{e}re symbol
satisfies property (\ref{eq:Steinberg}) - Theorem
\ref{th:Steinber-CC} - as an application of an adjunction formula
recently offered in \cite{Pa}, and, similar to other works by the
author, we offer a cohomological characterization of this symbol
from the commutator of a central extension of groups - Proposition
\ref{p:char}- in the same case.

\section{Preliminaries} \label{s:prel}

\subsection{Witt Parameters and the Lifting Lemma} \label{sub:Witt}

In this work {\em rings} are commutative with unit. Let $A$ be a
ring. Let $A((t))$ be the ring of series of the form $\sum_{i\in
\ZZ}a_it^i$ with coefficients $a_i\in A$ such that $a_i=0$ for
$i\ll 0$. Let $A[[t]]\subset A((t))$ be the subring consisting of
series of the form $\sum_{i=0}^\infty a_it^i$. Let $A[t^{\pm
1}]\subset A((t))$ be the subring consisting of polynomials in
$t^{\pm 1}$ with coefficients in $A$.  Given an ideal $I\subset
A$,
 let $t^{\pm 1}I[t^{\pm 1}]$
be the ideal of $A[t^{\pm 1}]$ generated by all products of the
form $xt^{\pm 1}$, where $x\in I$; let $I((t))$ be the ideal of
$A((t))$ consisting of series with all coefficients in $I$, and
let $I[[t]]=A[[t]]\cap I((t))$.

Let $A$ again be a ring. Let $I$ be a nilpotent ideal.  Let
$\Gamma(A,I)$ be the set of power series
$f=\sum_{i=-\infty}^\infty a_it^i\in A((t))$ such that for some
integer $w=w(f)=w_{A,I}(f)$ we have $a_w\in A^\times$ and $a_i\in
I$ for $i<w$. The set $\Gamma(A,I)$ is closed under power series
multiplication and forms a group. Let $\Gamma_0(A,I)\subset
\Gamma(A,I)$ be the subgroup consisting of power series for which
$w(f)=0$. Let $\Gamma_-(A,I)\subset \Gamma_0(A,I)$ be the subgroup
consisting of power series of the form $1+f$ with $f\in
t^{-1}I[t^{-1}]$. Let $\Gamma_+(A,I)=A[[t]]^\times\subset
\Gamma_0(A,I)$. Given pairs $(A,I)$ and $(B,J)$, each consisting
of a ring and a nilpotent ideal, and a ring homomorphism
$\varphi:A\rightarrow B$ such that $\varphi(I)\subset J$, the
corresponding group homomorphism
$\Gamma(\varphi):\Gamma(A,I)\rightarrow\Gamma(B,J)$ is defined to
be that sending $\sum_i a_it^i$ to $\sum_i \varphi(a_i)t^i$. Thus,
the construction $\Gamma$ becomes a functor, and similarly the
related constructions $\Gamma_{\pm}$ and $\Gamma_0$ become functors.

For a detailed study of the functor $\Gamma$, the reader is
referred to \cite{MPl}.

\begin{lem}
Let $A$ be a ring. Let $I\subset A$ be a nilpotent ideal. Let
$f\in \Gamma_0(A,I)$ be given. Then there exist unique $g\in
\Gamma_+(A,I)$ and $h\in \Gamma_-(A,I)$ such that $f=gh$.
\end{lem}
\begin{proof} Let $\nu$ be a positive integer such that $I^\nu=0$. We
proceed by induction on $\nu$. For $\nu=1$ there is nothing to
prove, so we assume that $\nu>1$ for the rest of the proof. Write
$f=f_++f_-$, where $f_+\in A[[t]]^\times$ and $f_-\in
t^{-1}I[t^{-1}]$, in the unique possible way. After replacing $f$
by $(f_+)^{-1}f$, we may assume without loss of generality that
$f=1-\tilde{f}$, where $\tilde{f}\in I((t))$. Write
$\tilde{f}=\tilde{f}_++\tilde{f}_-$, where $\tilde{f}_+\in I[[t]]$
and $\tilde{f}_-\in t^{-1}I[t^{-1}]$, in the only possible way.
After replacing $f$ by
$$\left(1+\sum_{i=1}^{\nu-1}(\tilde{f}_-)^i\right)(1-\tilde{f}_+)^{-1}f=(1-
\tilde{f}_-)^{-1}(1-\tilde{f}_+)^{-1}f,$$ we may assume without
loss of generality that $1-f\in I^2((t))$, in which case we are
finished by induction on $\nu$. \end{proof}

\begin{lem}
Let $A$ be a ring. Let $f\in A[[t]]^\times$ be given. Then:

(i) There exists a unique sequence $\{a_i\}_{i=1}^\infty$ in $A$
such that \linebreak $f=f(0)\prod_{i=1}^\infty (1-a_it^i)$.

(ii) If $f=1+t^{n+1}g$ for some $g\in A[[t]]$ and positive integer
$n$, then $a_1=\cdots=a_n=0$.

(iii) If there exists a nilpotent ideal $I\subset A$ such that
$1-f\in tI[t]$, then $a_i\in I$ for all $i$ and $a_i=0$ for all
$i\gg 0$.
\end{lem}

\begin{proof} (i) Write $f/f(0)=1-tg_1$ for some $g_1\in A[[t]]$.
 For $i>0$ there exists unique $g_i\in A[[t]]$
 such that
$$1-t^{i}g_i=(1-g_{i-1}(0)t^{i-1})^{-1}(1-t^{i-1}g_{i-1}).$$
Put $a_i=g_{i}(0)$ for all $i$. The resulting sequence
$\{a_i\}_{i=1}^\infty$ is the only possible one with the desired
properties.

(ii) The method of construction of the sequence
$\{a_i\}_{i=1}^\infty$ also proves this.

(iii) Let $\nu$ be a positive integer such that $I^\nu=0$. If
$\nu=1$, there is nothing to prove, so we assume that $\nu>1$ for
the rest of the proof. Write $f=1+\sum_{i=1}^\mu b_it^i$ for some
coefficients $b_i\in I$ and positive integer $\mu$.
 By the uniqueness asserted in part (i) we have congruences $a_i\equiv b_i\mod I^2$ for $i=1,\dots,\mu$
and $a_i\equiv 0\bmod I^2$ for $i>\mu$. Consider $f^*=
f\prod_{i=1}^\mu(1+\sum_{j=1}^{\nu-1}a_i^\nu
t^{i\nu})=f\prod_{i=1}^\mu(1-a_it^i)^{-1}$. Then,
$f^*=1-\sum_{\mu<i\leq \mu^*}b^*_it^i$ for some integer
$\mu^*\geq\mu$ and coefficients $b^*_i\in I^2$. Writing
$f^*=\prod_{i=1}^\infty (1-a^*_it^i)$ according to part (i) for
unique coefficients $a^*_i\in A$, we have $a^*_i=0$ for $i\gg 0$
by induction on $\nu$. We further have $a^*_1=\cdots=a^*_\mu=0$ by
part (ii), and finally we have $a^*_i=a_i$ for $i>\mu$ by the
uniqueness asserted in part (i). \end{proof}

\begin{prop}\label{p:Witt}
Let $A$ be a ring. Let $I\subset A$ be a nilpotent ideal. Let
$f\in \Gamma(A,I)$ be given. Then there exist unique coefficients
$\{a_i\}_{i=-\infty}^\infty$ in $A$ such that $a_0\in A^\times$,
$a_i\in I$ for $i<0$, $a_i=0$ for $i\ll 0$, and
$f=t^{w(f)}a_0(\prod_{i=1}^\infty(1-a_it^i))(\prod_{i=1}^\infty(1-a_{-i}t^{-i}))$.
\end{prop}

\noindent We call $\{a_i\}_{i=-\infty}^\infty$ the family of {\em
Witt parameters} $f\in \Gamma(A,I)$.

\begin{proof} We combine the
preceding two lemmas. \end{proof}

    We should note that the Witt parameters $\{a_i\}_{i=-\infty}^\infty$
depend functorially on $f$, i.~e., given a ring homomorphism
$\varphi:A\rightarrow A'$ and nilpotent ideals $I\subset A$ and
$I'\subset A'$ such that $\varphi(I)\subset I'$, if
$\{a_i\}_{i=-\infty}^\infty$ are the Witt parameters in $A$ of
$f\in \Gamma(A,I)$, then $\{\varphi(a_i)\}_{i=-\infty}^\infty$ are
the Witt parameters in $A'$ of $\Gamma(\varphi)(f)\in
\Gamma(A',I')$.

\begin{lem}[Lifting Lemma]\label{l:lift}
Let $A$ be a ring and let $I\subset A$ be a nilpotent ideal. Let
$f\in A((t))$ be given such that $f,1-f\in \Gamma_0(A,I)$. We can
then find:  a ring $A_1$ and nilpotent ideal $I_1\subset A_1$; an
artinian local $\QQ$-algebra $A_2$ and nilpotent ideal $I_2\subset
A_2$; a ring homomorphism $\varphi:A_1\rightarrow A$ such that
$\varphi(I_1)\subset I$;  an element $f_1\in A_1((t))$ such that
$f_1,1-f_1\in \Gamma_0(A_1,I_1)$, $\Gamma(\varphi)(f_1)=f$, and
$\Gamma(\varphi)(1-f_1)=1-f$, and a ring homomorphism \linebreak
$\psi:A_1\rightarrow A_2$ such that $\psi(I_1)\subset I_2$ and
$\Gamma(\psi)$ is injective.
\end{lem}

\begin{proof} Let us write $f=\sum_{i\in \ZZ}b_it^i$ where $b_i\in A$, $b_i\in
I$ for $i<0$, $b_i=0$ for $i\ll 0$. Let us choose a positive
integer $N$ such that $f=\sum_{i\geq -N}b_it^i$ and $b_i^N=0$ for
$i<0$. Let $\Lambda_0$ be the polynomial ring over $\ZZ$ generated
by a family of independent variables $\{X_i\}_{i\in \ZZ}$, and let
$\Lambda$ be obtained from $\Lambda_0$ by inverting $X_0$ and
$1-X_0$. Let $J_0\subset \Lambda$ be the ideal generated by
$$\{X_i\}_{i<-N}\cup\{X_i^N\}_{-N\leq i\leq -1}$$
and let $J_1\subset \Lambda$ be the ideal generated by
$\{X_i\}_{i<0}$. Let $A_1=\Lambda/J$ and $I_1=J_0/J_1$. Let $A_2$
be the ring obtained from $A_1$ by inverting all non-zero-divisors
(the total quotient ring of $A_1$) and let $I_2$ be the ideal of
$A_2$ generated by fractions with numerator in $I_1$. Let
$\varphi$ be the unique homomorphism sending $X_i$ to $a_i$ for
all $i$. Let $f_1=\sum_i X_it^i\in \Gamma(A_1,I_1)$. Let
$\psi:A_1\rightarrow A_2$ be the natural inclusion. Then the
objects thus constructed have all desired properties. \end{proof}

\subsection{Commensurability of $A$-modules and the Contou-Carr\`{e}re
symbol}\label{ss:central-extension}
\smallskip
Let $A$ now be an artinian local ring with maximal ideal $m$. Let
$V_A$ be a free $A$-module. Given free $A$-submodules $E,F\subset
V$, we write $E\sim F$ and say that $E$ and $F$ are {\em
commensurable} if the quotient $\frac{E+F}{E\cap F}$ is finitely
generated over $A$. It is easily verified that commensurability is
an equivalence relation.

    Given a free $A$-submodule $V_A^+\subset
V_A$, let $\glva$ denote the set of $A$-linear automorphisms
$\sigma$ of $V_A$ such that $V_A^+\sim \sigma V_A^+$. It is easily
verified that $\glva$ is a subgroup of the group of $A$-linear
automorphisms of $V_A$ depending only on the commensurability
class of $V_A^+$.

    Using the theory of groupoids, we can construct a group
 $\glvae$ depending only on the commensurability class of
 $V_A^+$ and such that there exists a central extension of
 groups:
 \begin{equation} \label{e:greg} 1 \to A^\times \to \glvae \overset \pi \longrightarrow \glva \to 1\,
 .\end{equation}

    We denote by $\{\cdot, \cdot\}_{V_A}^{V_A^+}$ the
 commutator of the central extension (\ref{e:greg}); that is, if $\tau$ and $\sigma$ are two
commuting elements of $\glva$ and ${\tilde \tau}, {\tilde
\sigma}\in \glvae$ are elements such that $\pi ({\tilde \tau}) =
\tau$ and $\pi ({\tilde \sigma}) = \sigma$, then one has a
commutator pairing:
$$\{\tau,\sigma\}_{V_A}^{V_A^+} = {\tilde \tau}\cdot {\tilde \sigma}\cdot {\tilde \tau}^{-1} \cdot {\tilde
\sigma}^{-1} \in A^\times\, .$$

    If we set $V_A^+ = A[[t]]$ and $V_A = A((t))$, we have that $A((t))^\times \subseteq
 \glva$ (viewing a Laurent series $f$ as a homothety) and $A((t))^\times$ is a commutative group.
From this immersion of groups, the central extension
(\ref{e:greg}) induces a new central extension of groups:
 \begin{equation} \label{e:est-fer} 1 \to A^\times \to {\widetilde {A((t))^\times}} \to A((t))^\times \to 1\,
 ,\end{equation}\noindent and we have a commutator map:
$$\{\cdot, \cdot\}_{A((t))}^{A[[t]]} \colon A((t))^\times \times
A((t))^\times\longrightarrow A^\times\, .$$

    Given $f\in A((t))^\times$, by considering the family of {\em Witt
parameters} of $f$ (Proposition \ref{p:Witt}), one has a unique
presentation:
$$f=a_0\cdot t^{w(f)}\cdot \prod_{i=1}^\infty
\left(1-a_{-i}t^{-i}\right) \cdot \prod_{i=1}^\infty
\left(1-a_it^i\right)\, ,
$$
where
$$w(f)\in \ZZ,\;\;\;
\left\{\begin{array}{ll}
a_i=0&\mbox{if $i\ll 0$,}\\a_i\in m&\mbox{if $i<0$,}\\
a_i\in A^\times&\mbox{if $i=0$,}\\
a_i\in A&\mbox{if $i>0$.}
\end{array}\right.$$
 The integer number $w(f)$ is the {\em winding number} of $f$.

    Thus, if $g\in A((t))^\times$ is another element with
presentation
$$g = b_0\cdot
t^{w(g)} \cdot {\prod}_{i=1}^{\infty} (1-b_{i}t^{i}) \cdot
{\prod}_{i=1}^{\infty} (1-b_{-i}t^{-i})\, ,$$\noindent such that
the coefficients satisfy the above conditions, a computation shows
that the value of the commutator is:
$$\{f, g\}_{A((t))}^{A[[t]]} = \frac
{a_0^{w(g)}{\prod}_{i=1}^{\infty}{\prod}_{j=1}^{\infty}\big (1 -
a_i^{j/(i,j)}b_{-j}^{i/(i,j)}\big
)^{(i,j)}}{b_0^{w(f)}{\prod}_{i=1}^{\infty}{\prod}_{j=1}^{\infty}\big
(1 - a_{-i}^{j/(i,j)}b_{j}^{i/(i,j)}\big )^{(i,j)}} \in A^\times\,
.$$

    This expression makes sense because only finitely many of the terms appearing in
the infinite products differ from 1, and  the Contou-Carr\`{e}re
symbol \cite{Co} is: $$\langle f,g\rangle_{A((t))^\times} =
(-1)^{w(f)w(g)}\cdot \{f, g\}_{A((t))}^{A[[t]]}\, .$$

    For arbitrary elements $f,g,h \in A((t))^\times$, the following relations hold:
\begin{itemize}

\item $\langle f,g\cdot h\rangle_{A((t))^\times} = \langle
f,g\rangle_{A((t))^\times} \cdot \langle f,
h\rangle_{A((t))^\times}$.

\item $\langle g,f\rangle_{A((t))^\times} = \langle
f,g\rangle_{A((t))^\times}^{-1}$.

\item $\langle f,-f\rangle_{A((t))^\times} = 1$.

\item Given $\varphi \in A((t))^\times$ with positive winding
number $n$, one has that: \begin{equation}\label{eq:norm} \langle
f,g \circ \varphi\rangle_{A((t))^\times} = \langle{\mathcal
N}_{\varphi} [f],g\rangle_{A((t))^\times}\,
,\end{equation}\noindent
 where ${\mathcal
N}_{\varphi}\colon A((t))^\times \to A((t))^\times$ denotes the
corresponding norm mapping -viewing $A((t))$ via the homomorphism
$h \mapsto h\circ \varphi$ as a free $A((t))$-module of rank $n$-
(\cite{Pa}, Proposition 3.6). As an application of this
``adjunction formula'' one has that the Contou-Carr\`{e}re symbol
is invariant under reparameterization of $A((t))$ in the following
sense: if $\tau\in A((t))$ is a element with winding number equal
to $1$, then $\langle f,g\rangle_{A((t))^\times} = \langle f\circ
\tau,g\circ \tau \rangle_{A((t))^\times}$.
\end{itemize}

\section{The Steinberg property of the Contou-Carr\`{e}re symbol}

    With the above notation, let us now consider a Laurent series $$f = \sum_{i\geq -N} a_i t^i \in A((t))^{\times}\,
.$$

\begin{lem} \label{l:inv} If $w(f) = w(1-f) = 0$ and $a_i \in A^\times$ for a
positive integer \linebreak $i > 0$, then there exists an
invertible element $\lambda \in A^\times$ and a series $\varphi
\in A((t))^\times$ with $w(\varphi) > 0$ such that $$1-f =
(1-\lambda)(1-\varphi)\, .$$
\end{lem}

\begin{proof} With the conditions of the Lemma, it is clear that
there exists an invertible element $\lambda \in A^\times$ and a
series $g \in A((t))^\times$ with $w(g) > 0$ such that $f =
\lambda \cdot (1+g)$.

    Thus, $$1 - f = (1 - \lambda) [1 - \frac{\lambda g}{1 - \lambda}]\, ,$$\noindent and writing
$\varphi = \frac{\lambda g}{1 - \lambda}$ the claim is deduced.
\end{proof}

    With the notation of the preceding Lemma, note that $f = \lambda +
(1-\lambda)\varphi$.

\begin{lem}\label{l:second} If $A$ is an artinian local $\QQ$-algebra, and ${\tilde f} \in A((t))$ is a
nilpotent element, one has that $$\langle 1 + {\tilde f}, 1 + \mu
{\tilde f} \rangle_{A((t))^\times} = 1$$\noindent for all $\mu \in
A^\times$.
\end{lem}

\begin{proof} Recall from (\cite{De}, p. 154) that if $A$ is a
$\QQ$-algebra and $g\in 1+m((t))$, then
$$\langle f,g\rangle_{A((t))^\times} = \exp(\Res_{t=0} [\log g\cdot d\log f])$$\noindent for all $g\in A((t))^\times$.

    Thus, the statement of the Lemma follows from the well-known
property of residues: $$\Res_{t=0} [{\tilde f}^n d{\tilde f}] =
0\, ,$$\noindent for every ${\tilde f} \in A((t))$ and $n\geq 0$.
\end{proof}

\begin{thm}\label{th:Steinber-CC} [Steinberg property] If $A$ is an artinian local ring,
given an element $f\in A((t))^\times$ such that $1-f \in
A((t))^\times$, one has that:
$$\langle f,1-f\rangle_{A((t))^\times} = 1\, .$$
\end{thm}

\begin{proof} As in the proof of J. Milnor \cite{Mi} related to the tame
symbol, the proof of this theorem will be divided into several
cases.

    If $w(f) = n > 0$, from expression
(\ref{eq:norm}) one has that:
$$\begin{aligned} \langle 1-f,f\rangle_{A((t))^\times} &= \langle 1-f,t\circ
f\rangle_{A((t))^\times}
\\ &= \langle {\mathcal N}_{f}
[1-f],t\rangle_{A((t))^\times} \\ &= \langle
(1-t)^n,t\rangle_{A((t))^\times} = 1 \, ,\end{aligned}$$\noindent
and the claim is deduced in this case.

    Moreover, when $w(f) < 0$, bearing in
mind that $\langle f^{-1},-f^{-1}\rangle_{A((t))^\times} = 1$,
from the above result we have that:
$$\begin{aligned} \langle f,1-f\rangle_{A((t))^\times} &= \langle
f^{-1},1-f\rangle_{A((t))^\times}^{-1}\cdot \langle
f^{-1},-f^{-1}\rangle_{A((t))^\times}^{-1}
\\ &=  \langle
f^{-1},(1-f)(-f^{-1})\rangle_{A((t))^\times}^{-1} \\ &= \langle
f^{-1}, 1-f^{-1}\rangle_{A((t))^\times}^{-1} = 1 \,
.\end{aligned}$$

    Furthermore, if $w(f) = w(1-f) = 0$, and $f$ satisfies
the condition of Lemma \ref{l:inv}, with the notation of this
Lemma one has that:
$$\begin{aligned} \langle f,1-f\rangle_{A((t))^\times} &= \langle
f,1-\lambda\rangle_{A((t))^\times} \cdot \langle f,1 - \varphi
\rangle_{A((t))^\times}
\\ &=   \langle f, (1 - t) \circ
\varphi \rangle_{A((t))^\times} = \langle {\mathcal N}_{\varphi}
[f], 1 - t \rangle_{A((t))^\times} \\ &= \langle [\lambda +
(1-\lambda)t]^{w(\varphi)}, 1 - t \rangle_{A((t))^\times} =  1 \,
.\end{aligned}$$

    Finally, if $w(f) = w(1-f) = 0$, and $a_i \in m$ for
all $i \ne 0$, it follows from the Lifting Lemma (Lemma
\ref{l:lift}) that we can assume without loss of generality that
$A$ is an artinian local $\QQ$-algebra.  Then, writing $f = a_0(1
+ {\tilde f})$ we have that $1 - f = (1-a_0)(1 -
\frac{a_0}{1-a_0}{\tilde f})$.

    Hence, setting $\mu = - \frac{a_0}{1-a_0}$, we conclude from Lemma \ref{l:second} bearing in mind that:
$$\langle f, 1-f \rangle_{A((t))^\times} = \langle 1 + {\tilde f}, 1 + \mu {\tilde f}
\rangle_{A((t))^\times} = 1\, .$$
\end{proof}

\section{Cohomological characterization of the Contou-Carr\`{e}re symbol}

    Let $\{\cdot, \cdot\}_{A((t))}^{A[[t]]}$ again be the
commutator of the central extension of groups (\ref{e:est-fer}).
Since $\{\cdot, \cdot\}_{A((t))}^{A[[t]]}$ is a 2-cocycle, it
determines an element of the cohomology group $H^2(A((t))^\times,
A^\times)$ -\cite{Se}, page 168-.

    Similar to above, we can say that a map $\psi \colon A((t))^{\times} \times
A((t))^{\times}\longrightarrow A^{\times}$ is called a ``Steinberg
map'' when:
\begin{itemize}
\item $\psi$ is bimultiplicative.
\item $\psi (f,1-f) = 1$ for all
$f\in A((t))^\times$ such that $1-f \in A((t))^\times$.
\end{itemize}

\begin{rem} A Steinberg map $\psi$ also satisfies the following
properties:
\begin{itemize}
\item $\psi (f,g) = [\psi (g,f)]^{-1}$; \item $\psi (f,-f) = 1$,
\end{itemize}
for all $f,g\in A((t))^\times$.
\end{rem}

\begin{rem} The commutator $\{\cdot, \cdot\}_{A((t))}^{A[[t]]}$ is
not a Steinberg map because $$\{\frac1t, 1-
\frac1t\}_{A((t))}^{A[[t]]}= -1\, .$$
\end{rem}

    We shall now give a cohomological characterization of the Contou-Carr\`{e}re
symbol from the cohomology class $[\{\cdot,
\cdot\}_{A((t))}^{A[[t]]}] \in H^2(A((t))^\times, A^\times)$.

    Recall that given two commutative groups, $G$ and $B$,
$H^2(G,B)$ is the group of classes of
 2-cocycles $c\colon G\times G \to B$ (mod. 2-coboundaries), where
a 2-cocycle
 $b$ is a 2-coboundary when there exists a map
 $\phi \colon  G \to B$
such that
 $$b(\alpha, \beta) = (\delta \phi) (\alpha,\beta) = \phi (\alpha \cdot \beta) \cdot {\phi (\alpha)}^{-1}
\cdot {\phi (\beta)}^{-1}\, .$$

Similar to \cite{Pa1}, one has that:

\begin{lem}\label{l:cobor} There exists a unique 2-coboundary
$$b \colon {\mathbb Z}
\times {\mathbb Z} \to A^{\times}$$ \noindent satisfying the
conditions:
\begin{itemize}
\item  $b(\alpha, \beta + \gamma) = b(\alpha, \beta)\cdot
 b(\alpha, \gamma)$;
\item $b(\alpha, \alpha) =  (-1)^{\alpha}$,
\end{itemize}
\noindent for all $\alpha, \beta, \gamma \in {\mathbb Z}$.
\end{lem}

\begin{proof} Let $\phi (\alpha)
= {\lambda}_{\alpha} \in A^{\times}$. If $b = \delta \phi$, it follows
from the above conditions that
$${\lambda}_{\alpha} = (-1)^{\frac {\alpha (\alpha-1)}{2}}
{\lambda}_{1}^{\alpha} \qquad \text { for each } \alpha\in
{\mathbb Z}\, ,$$\noindent and $b(\alpha, \beta) = (-1)^{\alpha
\cdot \beta}$ is therefore the unique 2-coboundary that satisfies
the required conditions.
\end{proof}

\begin{prop}[Cohomological characterization of the Contou-Carr\`{e}re symbol] \label{p:char} If $A$ is an
artinian local ring, there exists a unique Steinberg map $\langle
\quad,\quad \rangle_{A((t))^\times}$ in the cohomology class
 $$[\{\cdot,\cdot\}_{A((t))}^{A[[t]]}] \in H^2(A((t))^\times,
A^\times)$$\noindent satisfying the condition:
 $$\langle f,g \rangle_{A((t))^\times} = \{f,g\}_{A((t))}^{A[[t]]} \text{ if } w(f) = 0$$
\noindent for all $f,g \in A((t))^\times$.

This element is the Contou-Carr\`{e}re symbol $\langle \quad,\quad
\rangle_{A((t))^\times}$.
\end{prop}

\begin{proof}  Let $\nu (f,g) = c'(f,g)\cdot
\{f,g\}_{A((t))}^{A[[t]]}$ be an element of the cohomology class
$[\{\cdot,\cdot\}_{A((t))}^{A[[t]]}] \in H^2(A((t))^\times,
A^\times)$ satisfying the condition of the proposition. Since $c'$
is a 2-coboundary, one has that $c'(f,g)=1$ when $w(g)=0$ and,
therefore, there exists a commutative diagram
$$\xymatrix{ A((t))^{\times}\times  A((t))^{\times}
\ar[d]_{w\times w} \ar[dr]^{c'}  \\  {{\mathbb Z}\times {\mathbb
Z}} \ar[r]^{{\tilde c}'} & A^{\times}}$$ \noindent where ${{\tilde
c}'}$ is a 2-coboundary satisfying
$${{\tilde c}'}(x, y + z) = {{\tilde c}'}(x,
y)\cdot {{\tilde c}'}(x, z).$$

    Furthermore, since $\nu$ is a Steinberg map, then $\nu(f,-f)=1$, and one has that ${{\tilde c}'}(x, x) =
(-1)^{x}$.

    It then follows from Lemma \ref{l:cobor} that ${{\tilde c}'}(x,
y)  = (-1)^{x\cdot y}$ and $$c'(f,g) = (-1)^{w (f)\cdot w (g)}\,
.$$

Thus, the unique element in $[\{\cdot,\cdot\}_{A((t))}^{A[[t]]}]
\in H^2(A((t))^\times, A^\times)$ satisfying the condition of the
proposition is $\nu (f,g) = (-1)^{w(f) \cdot w(g)} \cdot
\{f,g\}_{A((t))}^{A[[t]]}$ and $\nu (f,g) = \langle f,g
\rangle_{A((t))^\times}$.
\end{proof}

\begin{rem} The above property, which characterizes the Contou-Carr\`{e}re symbol, is
equivalent to one of the conditions that J.P. Serre (\cite{Se})
gave to define local symbols on algebraic curves that are also
Steinberg symbols.

    In the formalism of the cohomological characterization of other
classical symbols, the author has formulated a conjecture to
replace this condition with the continuity of the map by
considering natural topologies on the corresponding groups $G$ and
$B$ that determine the cohomology group $H^2(G,B)$. However, in
this case we are not sufficiently confident to formulate a
conjecture about a finer characterization of the
Contou-Carr\`{e}re symbol.
\end{rem}
\bigskip

{\centerline {\bf ACKNOWLEDGMENTS}}

\medskip

The author thanks Prof. Greg W. Anderson for helpful personal
communication on Witt parameters and the lifting lemma that
contains the results offered in Subsection \ref{sub:Witt}.

\end{document}